\documentclass[12pt]{article}
\usepackage{amsfonts}

\usepackage{mathrsfs}
\usepackage{amscd}
\usepackage{amsmath,amsfonts,amssymb,amscd}
\usepackage{indentfirst,graphics,epsfig,psfrag}
\input{epsf}

\usepackage{amsmath,amssymb,amsthm, amscd, amsfonts, graphicx,ifpdf}
\usepackage{lineno}

\newtheorem{thm}{Theorem}[section]

\newtheorem{lem}{Lemma}[section]

\def\qed{\nopagebreak\hfill{\rule{4pt}{7pt}}
\medbreak}

\def\pf{\noindent {\it Proof.} }

\title{\bf Note on the minimal size of a graph\\ with generalized connectivity
$\kappa_{3}= 2$\footnote{Supported by NSFC and the Fundamental
Research Funds for the Central Universities. }}

\author{
\small Shasha Li, Xueliang Li, Yongtang Shi\\
\small Center for Combinatorics and LPMC-TJKLC\\
\small Nankai University, Tianjin 300071, China.\\
\small  Email: lss@cfc.nankai.edu.cn, lxl@nankai.edu.cn, shi@nankai.edu.cn\\
}
\date{}
\begin{document}

\maketitle

\begin{abstract}

The concept of generalized $k$-connectivity $\kappa_{k}(G)$ of a
graph $G$ was introduced by Chartrand et al. in recent years. In our
early paper, extremal theory for this graph parameter was started.
We determined the minimal number of edges of a graph of order $n$
with $\kappa_{3}= 2$, i.e., for a graph $G$ of order $n$ and size
$e(G)$ with $\kappa_{3}(G)= 2$, we proved that $e(G)\geq
\frac{6}{5}n$, and the lower bound is sharp by constructing a class
of graphs, only for $n\equiv 0 \ (mod \ 5)$ and $n\neq 10$. In this
paper, we improve the lower bound to $\lceil\frac{6}{5}n\rceil$.
Moreover, we show that for all $n\geq 4$ but $n= 9, 10$, there
always exists a graph of order $n$ with $\kappa_{3}= 2$ whose size
attains the lower bound $\lceil\frac{6}{5}n\rceil$. Whereas for $n=
9, 10$ we give examples to show that $\lceil\frac{6}{5}n\rceil+1$ is
the best possible lower bound. This gives a clear picture on the
minimal size of a graph of order $n$ with generalized connectivity
$\kappa_{3}= 2$.
\\[3mm]
\end{abstract}

\section{Introduction}

We follow the terminology and notations of \cite{Bondy}, and all
graphs considered here are always finite and simple. As usual, we
denote the numbers of vertices and edges in $G$ by $n(G)$ and $e(G)$
(or simply $n$ and $e$), and these two basic parameters are called
the $order$ and $size$ of $G$, respectively. A stable set in a graph
is a set of vertices no two of which are adjacent. A vertex with
degree one in a tree is called a leaf. The $connectivity$
$\kappa(G)$ of a graph $G$ is defined as the minimum cardinality of
a set $Q$ of vertices of $G$ such that $G-Q$ is disconnected or
trivial. A well-known theorem of Whitney \cite{Whitney} provides an
equivalent definition of the connectivity. For each $2$-subset
$S=\{u,v\}$ of vertices of $G$, let $\kappa(S)$ denote the maximum
number of internally disjoint $uv$-paths in $G$. Then
$\kappa(G)=$min$\{\kappa(S)\}$, where the minimum is taken over all
$2$-subsets $S$ of $V(G)$.

In \cite{Chartrand}, the authors generalized the concept of
connectivity as follows. Let $G$ be a nontrivial connected graph of
order $n$ and $k$ an integer with $2\leq k\leq n$. For a set $S$ of
$k$ vertices of $G$, let $\kappa (S)$ denote the maximum number
$\ell$ of edge-disjoint trees $T_1,T_2,\ldots,T_\ell$ in $G$ such
that $V(T_i)\cap V(T_j)=S$ for every pair $i,j$ of distinct integers
with $1\leq i,j\leq \ell$ (note that the trees are vertex-disjoint
in $G\backslash S$). The $k$-$connectivity$, denoted by
$\kappa_k(G)$, of $G$ is then defined by
$\kappa_k(G)=$min$\{\kappa(S)\}$, where the minimum is taken over
all $k$-subsets $S$ of $V(G)$. Obviously, $\kappa_2(G)=\kappa(G)$.

This paper is a further development of our early work \cite{LLS},
where we determined the minimal number of edges of a graph with
$\kappa_{3}= 2$, i.e., for a graph $G$ of order $n$ and size $e(G)$
with $\kappa_{3}(G)= 2$, we proved that $e(G)\geq \frac{6}{5}n$, and
the lower bound is sharp by constructing a class of graphs, only
$n\equiv 0 \ (mod \ 5)$ and $n\neq 10$. Note that the number of
edges is integral and so the order of the graph attaining the lower
bound must be a multiple of $5$. On the other hand, since $e(G)$ is
an integer, the lower bound can be naturally improved to
$\lceil\frac{6}{5}n\rceil$. In this paper, we want to show that for
all $n\geq 4$ but $n= 9, 10$, the lower bound
$\lceil\frac{6}{5}n\rceil$ is best possible, whereas for $n= 9, 10$
we give examples to show that $\lceil\frac{6}{5}n\rceil+1$ is the
best possible lower bound. This gives a clear picture on the minimal
size of a graph of order $n$ with generalized connectivity
$\kappa_{3}= 2$.

\section{Preliminaries}

Before proceeding, we list some known results in \cite{LLS}
and \cite{LLZ}.

\begin{lem}[\cite{LLZ}]\label{lem1}
If $G$ is a connected graph with minimum degree $\delta$, then
$\kappa_3(G)\leq \delta$. In particular, if there are two adjacent
vertices of degree $\delta$, then $\kappa_3(G)\leq \delta-1$.
\end{lem}

\begin{lem}[\cite{LLS}]\label{lem2}
For a positive integer $k\neq 2$, let $C=x_{1}y_{1}x_{2}y_{2}\ldots
x_{2k}y_{2k}x_{1}$ be a cycle of length $4k$. Add $k$ new vertices
$z_{1},z_{2},\ldots,z_{k}$ to $C$, and join $z_{i}$ to $x_{i}$ and
$x_{i+k}$, for $1\leq i\leq k$. The resulting graph is denoted by
$H$. Then, the $3$-connectivity of $H$ is $2$, namely,
$\kappa_3(H)=2$.
\end{lem}

\begin{lem}[\cite{LLS}]\label{lem3}
For any connected graph $G$ of order $10$ and size $12$,
$\kappa_{3}(G)=1$.
\end{lem}

\noindent \textbf{Remark 2.1:} Note that there exists a graph $G$
such that $n=10$, $e(G)=13$ and $\kappa_{3}(G)=2$, as shown in
Figure $1$.
\begin{center}
\begin{picture}(150,90)
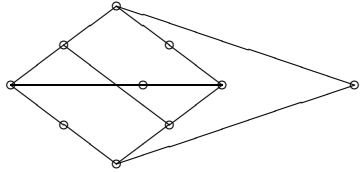
\put(-90,-10){Figure $1$: The graph $G$ of order $10$ and size $13$ with $\kappa_{3}(G)= 2$.}

\put(40,20){\circle{2.5}}
\put(20,35){\circle{2.5}}
\put(60,35){\circle{2.5}}
\put(0,50){\circle{2.5}}
\put(50,50){\circle{2.5}}
\put(80,50){\circle{2.5}}
\put(20,65){\circle{2.5}}
\put(40,80){\circle{2.5}}
\put(60,65){\circle{2.5}}
\put(130,50){\circle{2.5}}

\put(40,20){\line(4,3){20}}
\put(40,20){\line(-4,3){20}}
\put(60,35){\line(4,3){20}}
\put(20,35){\line(-4,3){20}}
\put(0,50){\line(4,3){20}}
\put(0,50){\line(1,0){50}}
\put(20,65){\line(4,3){20}}
\put(80,50){\line(-4,3){20}}
\put(60,65){\line(-4,3){20}}
\put(50,50){\line(1,0){30}}
\put(40,20){\line(3,1){90}}
\put(130,50){\line(-3,1){90}}
\put(60,35){\line(-4,3){40}}

\end{picture}
\end{center}

\vskip 3mm Now we turn to the graphs of order $9$ and size $11$.

\begin{lem}\label{lem4}
For any connected graph $G$ of order $9$ and size $11$, $\kappa_{3}(G)=1$.
\end{lem}
\pf Assume, to the contrary, that there is a connected graph $G$ of
order $n=9$ and size $m=11$ with $\kappa_{3}(G)=2$. By Lemma
\ref{lem1}, we have the minimum degree $\delta(G)\geq 2$. Denote by
$X$ the set of vertices of degree $2$ in $G$. It follows that
$2m=\Sigma_{v\in V(G)}d(v)\geq 2|X|+3(n-|X|)$, namely, $|X|\geq
3n-2m=5$. On the other hand, by Lemma \ref{lem1} again, we get that
$X$ is a stable set. Let $m'$ be the number of edges joining two
vertices belonging to $Y$, where $Y=V(G)-X$. It is clear that
$m=2|X|+m'$. So $|X|\leq \frac{m}{2}=5.5$. Now we can conclude that
$|X|=5$, $|Y|=4$, $m'=1$ and every vertex in $Y$ has degree exactly
$3$. Set $X=\{x_{1},x_{2},x_{3},x_{4},x_{5}\}$ and
$Y=\{y_{1},y_{2},y_{3},y_{4}\}$. Since $m'=1$, without loss of
generality, suppose that $y_1y_2$ is the only edge.

\noindent \textbf{Case 1:} There is a vertex in $X$ that is adjacent
to both $y_{1}$ and $y_{2}$.

Note that $G$ is a simple connected graph and every vertex in $X$
has degree $2$. It is not hard to get that $G$ is isomorphic to the
graph as shown in Figure $2$. Then observe that it is impossible to
find two internally-disjoint trees connecting the vertices $x_{1}$,
$x_{2}$ and $x_{4}$, contrary to our assumption.

\begin{center}
\begin{picture}(120,90)
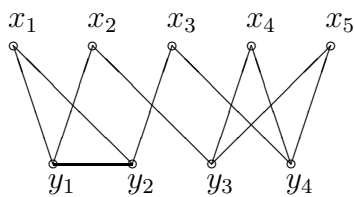
\put(-60,0){Figure $2$: The graph for Case $1$ of Lemma \ref{lem4}}

\put(13,21){$y_1$}
\put(43,21){$y_2$}
\put(73,21){$y_3$}
\put(103,21){$y_4$}

\put(-2,82){$x_1$}
\put(28,82){$x_2$}
\put(58,82){$x_3$}
\put(88,82){$x_4$}
\put(118,82){$x_5$}

\put(15,30){\circle{2.5}}
\put(45,30){\circle{2.5}}
\put(75,30){\circle{2.5}}
\put(105,30){\circle{2.5}}

\put(0,75){\circle{2.5}}
\put(30,75){\circle{2.5}}
\put(60,75){\circle{2.5}}
\put(90,75){\circle{2.5}}
\put(120,75){\circle{2.5}}

\put(15,30){\line(-1,3){15}}
\put(15,30){\line(1,3){15}}
\put(15,30){\line(1,0){30}}
\put(45,30){\line(-1,1){45}}
\put(45,30){\line(1,3){15}}
\put(75,30){\line(-1,1){45}}
\put(75,30){\line(1,3){15}}
\put(75,30){\line(1,1){45}}
\put(105,30){\line(-1,1){45}}
\put(105,30){\line(-1,3){15}}
\put(105,30){\line(1,3){15}}

\end{picture}
\end{center}

\noindent \textbf{Case 2:} There is no vertex in $X$ that is adjacent
to both $y_{1}$ and $y_{2}$.

\noindent \textbf{Subcase 2.1:} For every $2$-subset $\{y_{i},y_{j}\}$
of $Y$ other than $\{y_{1},y_{2}\}$, there is a vertex in $X$ that is
adjacent to both $y_{i}$ and $y_{j}$, where $1\leq i\neq j\leq 5$.

Note that there are exactly five vertices in $X$ and five
$2$-subsets of $Y$ other than $\{y_1,y_2\}$, namely,
$\{y_{1},y_{3}\},\{y_{1},y_{4}\},\{y_{2},y_{3}\},$
$\{y_{2},y_{4}\},\{y_{3},y_{4}\}$. Thus, we may assume that $G$ is
isomorphic to the graph as shown in Figure $3$. Consider the three
vertices $x_{1}$, $x_{2}$ and $x_{5}$, and we can get
$\kappa_3(G)=1$, contrary to our assumption.

\begin{center}
\begin{picture}(120,90)
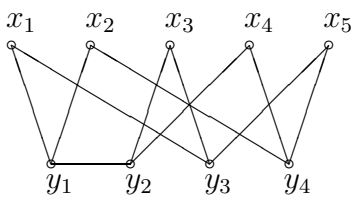
\put(-60,0){Figure $3$: The graph for Subcase $2.1$ of Lemma \ref{lem4}}

\put(13,21){$y_1$}
\put(43,21){$y_2$}
\put(73,21){$y_3$}
\put(103,21){$y_4$}

\put(-2,82){$x_1$}
\put(28,82){$x_2$}
\put(58,82){$x_3$}
\put(88,82){$x_4$}
\put(118,82){$x_5$}

\put(15,30){\circle{2.5}}
\put(45,30){\circle{2.5}}
\put(75,30){\circle{2.5}}
\put(105,30){\circle{2.5}}

\put(0,75){\circle{2.5}}
\put(30,75){\circle{2.5}}
\put(60,75){\circle{2.5}}
\put(90,75){\circle{2.5}}
\put(120,75){\circle{2.5}}

\put(15,30){\line(-1,3){15}}
\put(15,30){\line(1,3){15}}
\put(15,30){\line(1,0){30}}
\put(45,30){\line(1,3){15}}
\put(45,30){\line(1,1){45}}
\put(75,30){\line(-5,3){75}}
\put(75,30){\line(-1,3){15}}
\put(75,30){\line(1,1){45}}
\put(105,30){\line(-5,3){75}}
\put(105,30){\line(-1,3){15}}
\put(105,30){\line(1,3){15}}

\end{picture}
\end{center}

\noindent \textbf{Subcase 2.2:} Except $\{y_1,y_2\}$, there
exists another $2$-subset such that no vertex in $X$ is
adjacent to both of the vertices in that subset.

In such a situation, there must exist some $2$-subset
$\{y_{i},y_{j}\}$ such that at least two vertices in $X$ are
adjacent to both $y_{i}$ and $y_{j}$, where $1\leq i\neq j\leq 5$.
If $\{y_{i},y_{j}\}=\{y_{3},y_{4}\}$, it is not hard to get that
there must exist a vertex in $X$ that is adjacent to both $y_{1}$
and $y_{2}$, contrary to the case. So without loss of generality, we
may assume that $\{y_{i},y_{j}\}=\{y_{1},y_{3}\}$. Then we can get
$G$ is isomorphic to the graph as shown in Figure $4$. Observe that
it is impossible to find two internally-disjoint trees connecting
the vertices $x_{1}$, $x_{4}$ and $x_{5}$, contrary to our
assumption.

\begin{center}
\begin{picture}(120,90)
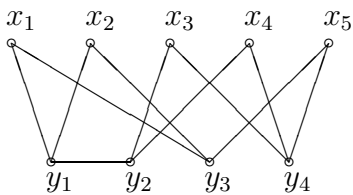
\put(-60,0){Figure $4$: The graph for Subcase $2.2$ of Lemma \ref{lem4}}

\put(13,21){$y_1$}
\put(43,21){$y_2$}
\put(73,21){$y_3$}
\put(103,21){$y_4$}

\put(-2,82){$x_1$}
\put(28,82){$x_2$}
\put(58,82){$x_3$}
\put(88,82){$x_4$}
\put(118,82){$x_5$}

\put(15,30){\circle{2.5}}
\put(45,30){\circle{2.5}}
\put(75,30){\circle{2.5}}
\put(105,30){\circle{2.5}}

\put(0,75){\circle{2.5}}
\put(30,75){\circle{2.5}}
\put(60,75){\circle{2.5}}
\put(90,75){\circle{2.5}}
\put(120,75){\circle{2.5}}

\put(15,30){\line(-1,3){15}}
\put(15,30){\line(1,3){15}}
\put(15,30){\line(1,0){30}}
\put(45,30){\line(1,3){15}}
\put(45,30){\line(1,1){45}}
\put(75,30){\line(-5,3){75}}
\put(75,30){\line(-1,1){45}}
\put(75,30){\line(1,1){45}}
\put(105,30){\line(-1,1){45}}
\put(105,30){\line(-1,3){15}}
\put(105,30){\line(1,3){15}}

\end{picture}
\end{center}
The proof is complete. \qed

\noindent \textbf{Remark 2.2:} Notice that there exists a graph $G$
such that $n=9$, $e(G)=12$ and $\kappa_{3}(G)=2$, as shown in Figure
$5$.
\begin{center}
\begin{picture}(150,90)
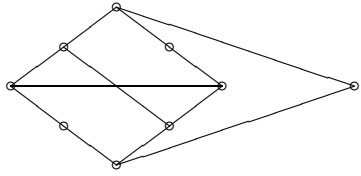
\put(-90,-10){Figure $5$: The graph $G$ of order $9$ and size $12$ with $\kappa_{3}(G)= 2$.}

\put(40,20){\circle{2.5}}
\put(20,35){\circle{2.5}}
\put(60,35){\circle{2.5}}
\put(0,50){\circle{2.5}}
\put(80,50){\circle{2.5}}
\put(20,65){\circle{2.5}}
\put(40,80){\circle{2.5}}
\put(60,65){\circle{2.5}}
\put(130,50){\circle{2.5}}

\put(40,20){\line(4,3){20}}
\put(40,20){\line(-4,3){20}}
\put(60,35){\line(4,3){20}}
\put(20,35){\line(-4,3){20}}
\put(0,50){\line(4,3){20}}
\put(0,50){\line(1,0){80}}
\put(20,65){\line(4,3){20}}
\put(80,50){\line(-4,3){20}}
\put(60,65){\line(-4,3){20}}
\put(40,20){\line(3,1){90}}
\put(130,50){\line(-3,1){90}}
\put(60,35){\line(-4,3){40}}

\end{picture}
\end{center}

\vskip 3mm Next we describe an operation on a vertex of degree $2$.

For a vertex $u$ of degree $2$, to {\it smooth $u$} is to delete $u$
and then add an edge between its neighbors. Obviously, performing
such an operation, the numbers of vertices and edges decrease by
one, respectively. Moreover, the degrees of the remaining vertices
are not changed.

\begin{lem}\label{lem5}
Let $G$ be a graph such that the set $X$ of vertices of degree $2$
is nonempty. Denote by $G'$ the new graph obtained by smoothing a
vertex in $X$, and then we have $\kappa_3(G')\geq \kappa_3(G)$.
\end{lem}
\pf Let $u$ be a vertex in $X$ and $\{w_1,w_2\}$ the neighbor set of
$u$. Suppose that $G'$ is obtained by smoothing $u$. Clearly,
$V(G')=V(G)-u$. For any three vertices $v_1$, $v_2$ and $v_3$ of
$G'$, let $S=\{v_1,v_2,v_3\}$. Obviously, $S\subseteq V(G)$. Let $T$
be a tree connecting $S$ in $G$. Note that if $v$ is a leaf of $T$,
we can assume that $v\in S$. Otherwise, $T'=T-v$ is still a tree
connecting $S$ and uses less vertices. Now if $u\in V(T)$, then we
can see that $T'=T-u+w_1w_2$ is exactly a tree connecting $S$ in
$G'$. If $u\notin V(T)$, the operation of smoothing $u$ has nothing
to do with $T$ and so $T$ is still a tree connecting $S$ in $G'$.
Therefore, it is not hard to get that $\kappa_{G'}(S)\geq
\kappa_{G}(S)$. From the definition of $\kappa_3$, the conclusion
that $\kappa_3(G')\geq \kappa_3(G)$ follows. \qed

\noindent \textbf{Remark 2.3:} For a given $G$, if we successively do the operation
of smoothing a vertex of degree $2$ more than once, the final graph
is denoted by $G'$. We can also get $\kappa_3(G')\geq \kappa_3(G)$.

\section{Lower bound}
\begin{lem}[\cite{LLS}]\label{lem6}
If $G$ is a graph of order $n$ with $\kappa_{3}(G)= 2$, then
$e(G)\geq \frac{6}{5}n$ and the lower bound is sharp.
\end{lem}

Note that the number of edges is integral and so the order of the
graph attaining the lower bound must be a multiple of $5$. In
\cite{LLS}, we showed that for all positive integer $k$ other than
$2$, there exists a graph of order $n=5k$ which attains the lower
bound. On the other hand, since $e(G)$ is an integer, the lower
bound can be improved to $\lceil\frac{6}{5}n\rceil$. Naturally, we
want to know whether there is a graph of order $n$ attaining the
lower bound for any positive integer $n$.

\begin{thm}\label{thm1}
If $G$ is a graph of order $n$ with $\kappa_{3}(G)= 2$, then
$e(G)\geq \lceil\frac{6}{5}n\rceil$. Moreover, the lower bound
is sharp for all $n\geq 4$ and $n\neq 9, 10$.
\end{thm}
\pf Since the number of edges must be an integer, by Lemma
\ref{lem6}, the lower bound $\lceil\frac{6}{5}n\rceil$ is obvious.

Note that all graphs considered here are always simple. Therefore, any graph
attaining the lower bound must have at least four vertices.
Moreover, by Lemmas \ref{lem3} and \ref{lem4}, we know that there is no
simple connected graph $G$ of order $9$ and size $11$ or
order $10$ and size $12$ such that $\kappa_3(G)=2$.

For $n=8$, there is a graph $G'$ of order $n$ such that
$\kappa_3(G')=2$ as shown in Figure $6$.  Moreover,
$e(G')=10=\lceil\frac{6}{5}\times 8\rceil$, which means that $G'$
attains the lower bound for $n=8$.

\begin{center}
\begin{picture}(150,90)
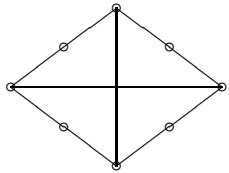
\put(-90,-10){Figure $6$: The graph $G'$ attaining the lower bound for $n=8$}

\put(40,20){\circle{2.5}}
\put(20,35){\circle{2.5}}
\put(60,35){\circle{2.5}}
\put(0,50){\circle{2.5}}
\put(80,50){\circle{2.5}}
\put(20,65){\circle{2.5}}
\put(40,80){\circle{2.5}}
\put(60,65){\circle{2.5}}

\put(40,20){\line(4,3){20}}
\put(40,20){\line(-4,3){20}}
\put(60,35){\line(4,3){20}}
\put(20,35){\line(-4,3){20}}
\put(0,50){\line(4,3){20}}
\put(0,50){\line(1,0){80}}
\put(20,65){\line(4,3){20}}
\put(80,50){\line(-4,3){20}}
\put(60,65){\line(-4,3){20}}
\put(40,20){\line(0,1){60}}

\end{picture}
\end{center}

\vskip 3mm Now, smooth a vertex of degree $2$ in $G'$. Clearly, the
resulting graph $G''$ is simple and $\delta(G'')=2$. By Lemma
\ref{lem5}, we can get $\kappa_3(G'')\geq (\kappa_(G')=2)$ and so
clearly $\kappa_3(G'')=2$. Moreover, $n=8-1=7$ and
$e=10-1=9=\lceil\frac{6}{5}\times 7\rceil$. The graph $G''$ is what
we want to find for $n=7$. Similarly, the graph obtained from $G''$
by smoothing any one vertex of degree $2$ attains the lower bound
for $n=6$.

Next, we consider the graph $H$ in Lemma \ref{lem2}. In \cite{LLS},
We obtained that $\kappa_3(H)=2$, $n(H)=5k$ and $e(H)=6k$, for
$k\neq 2$. So $H$ is exactly the graph of order $n=5k$ which attains
the lower bound.

For $k\geq 3$, let $k'=k-1$ and then $n(H)=5k'+5$ and $e(H)=6k'+6$.
Let $X$ be the set of vertices of degree $2$. Clearly $|X|=3k'+3>4$,
where $k'\geq 2$. Now for the graph $H$, smooth successively any $t$
vertices in $X$, for $1\leq t\leq 4$. For any $t$, it is easy to
check that no parallel edge can arise. Moreover, since $|X|>4$, the
minimum degree of the resulting graph $H'$ is still $2$. Combining
Lemma \ref{lem1} and Remark $2.3$, we can get the $3$-connectivity
of the resulting graph $H'$ is $2$. Now let us consider the numbers
of vertices and edges of $H'$.

When $t=1$, $n(H')=5k'+4$ and
$e(H')=6k'+5=\lceil\frac{6}{5}(5k'+4)\rceil$;

When $t=2$, $n(H')=5k'+3$ and
$e(H')=6k'+4=\lceil\frac{6}{5}(5k'+3)\rceil$;

When $t=3$, $n(H')=5k'+2$ and
$e(H')=6k'+3=\lceil\frac{6}{5}(5k'+2)\rceil$;

When $t=4$, $n(H')=5k'+1$ and
$e(H')=6k'+2=\lceil\frac{6}{5}(5k'+1)\rceil$.

Note that $k'\geq 2$. Therefore, for all $n\geq 4$ but $n= 9, 10$,
we can always find a graph of order $n$ attaining the lower bound.
\qed

\end{document}